\documentclass[12pt]{amsart}
\usepackage{amssymb}
\usepackage{amsbsy}
\usepackage{amscd}
\advance\textheight by \topskip
\usepackage{epic}
\usepackage{eepic}
\makeatletter
\usepackage{youngtab}
\def\cal{\mathcal}
\def\Bbb{\mathbb}
\def\frak{\mathfrak}

\newenvironment{pf*}[1]{\proof[#1]}{\endproof}
\renewcommand{\labelenumi}{(\theenumi)}%
\newcommand{\rom}{\textup}
\renewcommand{\thesubsection}{\thesection(\@roman\c@subsection)}
\makeatother
\newenvironment{aenume}{%
  \begin{enumerate}%
  \renewcommand{\theenumi}{\alph{enumi}}%
  \renewcommand{\labelenumi}{(\theenumi)}%
  }{\end{enumerate}}
\newtheorem{Theorem}[equation]{Theorem}
\newtheorem{Corollary}[equation]{Corollary}
\newtheorem{Lemma}[equation]{Lemma}
\newtheorem{Proposition}[equation]{Proposition}

\theoremstyle{definition}

\theoremstyle{remark}
\newtheorem{Remark}[equation]{Remark}

\numberwithin{equation}{section}
\numberwithin{figure}{section}

\newcommand{\thmref}[1]{Theorem~\ref{#1}}
\newcommand{\secref}[1]{\S\ref{#1}}
\newcommand{\lemref}[1]{Lemma~\ref{#1}}
\newcommand{\propref}[1]{Proposition~\ref{#1}}
\newcommand{\corref}[1]{Corollary~\ref{#1}}
\newcommand{\subsecref}[1]{\S\ref{#1}}

\newcommand{\defeq}{\overset{\operatorname{\scriptstyle def.}}{=}}
\newcommand{\C}{{\Bbb C}}
\newcommand{\Z}{{\Bbb Z}}
\newcommand{\Q}{{\Bbb Q}}
\newcommand{\R}{{\Bbb R}}

\newcommand{\CP}{\operatorname{\C P}}

\newcommand{\Ker}{\operatorname{Ker}}

\newcommand{\Ima}{\operatorname{Im}}

\newcommand{\pd}[2]{\frac{\partial#1}{\partial#2}}

\newcommand{\Hilb}[2]{{#1}^{\lbrack{#2}\rbrack}} 

\newcommand{\HilbX}[1]{\Hilb{X}{#1}}
\newcommand{\shfO}{\cal O} 
\newcommand{\idl}{\cal J} 
\newcommand{\Supp}{\operatorname{Supp}} 

\newcommand{\Jackl}{P_\lambda^{(\alpha)}}
\newcommand{\Jackm}{P_\mu^{(\alpha)}}

\newcommand{\barlambda}{{\bar\lambda}}
\newcommand{\barmu}{{\bar\mu}}
\newcommand{\HT}{H_{S^1}} 
\newcommand{\HTc}{H_{S^1,c}} 

\newcommand{\id}{\operatorname{id}}
\begin{document}
\title[Jack polynomials and Hilbert Schemes]
{Jack Polynomials and Hilbert Schemes of Points on Surfaces}
\author{Hiraku Nakajima}
\address{Department of Mathematical Sciences, University of Tokyo\\
Komaba 3-8-1, Meguro-ku, Tokyo 153, Japan}
\email{nakajima@ms.u-tokyo.ac.jp}
\thanks{Author supported in part by Grant-in-Aid for Scientific Research
(No.\ 08740010), Ministry of Education, Science and Culture, Japan.
}
\keywords{Jack polynomials, Hilbert scheme}
\subjclass{Primary 14C05, 05E05}
\maketitle

\section{Introduction}

The Jack (symmetric) polynomials $\Jackl(x)$ form a class of symmetric
polynomials which are indexed by a partition $\lambda$ and depend
rationally on a parameter $\alpha$. They reduced to the Schur
polynomials when $\alpha = 1$, and to other classical families of
symmetric polynomials for several specific parameters. Recently they
attracts attention from various points of view, for example the
integrable systems and combinatorics. They are simultaneous
eigenfunctions of certain commuting families of differential
operators, appearing in the integrable system called the
Calogero-Sutherland system (see e.g., \cite{AMOS} and the reference
therein). On the other hand, Macdonald studied their combinatorial
properties.  In fact, he introduced an even more general class of
symmetric functions (a two parameter family), which have many common
combinatorial features as Jack polynomials (see
\cite[Chapter~6]{Symmetric}).

It is well-known that Schur polynomials can be realized as certain
elements of homology groups of Grassmann manifolds (see e.g.,
\cite[Chapter~14]{Fulton}). The purpose of this paper is to give a
similar geometric realization for Jack polynomials. However, spaces
which we use are totally different. Our spaces are Hilbert schemes of
points on a surface $X$ which is the total space of a line bundle $L$
over the projective line $\CP^1$. The parameter $\alpha$ in Jack
polynomials relates to our surface $X$ by
\begin{equation}
  \alpha = -\langle{C},{C}\rangle = -c_1(L)[{C}],
\label{eq:selfinter}\end{equation}
where $C$ is the zero section, and $\langle{C},{C}\rangle$ is the
self-intersection number of $C$. It seems difficult to realize Jack
polynomials in homology groups of Grassmann manifolds since they
have no parameter.

The Hilbert scheme $\HilbX{n}$ parametrizing $0$-dimensional
subschemes of length $n$ on $X$ has been studied by various peoples
(see \cite{Go-book} and the reference therein). The author and
Grojnowski independently showed that the direct sum of the homology
groups of $\HilbX{n}$ (the summation is over $n$) is a representation
space of the Heisenberg algebra (boson Fock space) \cite{Gr,Na-hilb}.
It is also well-known that the space of symmetric polynomials is the
same representation space. This is the relationship between symmetric
functions and Hilbert schemes.

Now we explain our realization in more detail.
There are two characterizations of Jack polynomials:
\begin{aenume}
\item an orthogonal basis such that the transition matrix to monomial
symmetric functions is strictly upper triangular (see
\thmref{thm:defJack}), or
\item simultaneous eigenfunctions of a family of commuting
  differential operators (see above).
\end{aenume}
We use the characterization~(a). We first identify the complexified
ring of symmetric functions with the direct sum of the middle degree
homology groups of Hilbert schemes as above. We then identify the
inner product with the intersection pairing (\thmref{thm:ident}). This
result was essentially proved in
\cite{Gr,Na-hilb} combined with \cite{ES2}.
Then monomial symmetric functions are identified with fundamental
classes of certain middle dimensional subvarieties in the Hilbert
schemes (\thmref{thm:monom}), which were first introduced by
Grojnowski~\cite{Gr}. 
Finally, in order to get mutually orthogonal elements, we use the {\it
  localization\/}.  We define an $S^1$-action on the surface $X$ which
induces an action on the Hilbert scheme.  The localization, which goes
back to Bott residue formula \cite{Bott}, enables computations of the
intersection product to be reduced to the fixed point set of the
$S^1$-action.  Then cohomology classes which localize to different
fixed point components are mutually orthogonal. This is the mechanism
to construct a orthogonal basis. (In this paper, we shall use the
equivariant cohomology to formulate
the localization following \cite{AtBo}.)
The relation to the subvarieties corresponding to monomial symmetric
functions can be studied as follows. For each fixed point component,
we associate a locally closed submanifold (``stratum'') consisting of
points which converge to the fixed point components when they moved by
the $\C^*$-actions which extends the $S^1$-action (see
\eqref{eq:limit}). The dominance order is identified with the closure
relation of the strata \eqref{eq:closure}. The subvarieties
corresponding monomial symmetric functions are the closures of
stratum. Then it is easy to check that the transition matrix is
strictly upper triangular.

In fact, almost all arguments work for the total space of a line
bundle over any compact Riemann surface. We use the middle degree
homology group of Hilbert schemes on which only $H^2(C)$
contributes. ($H^0(C)$ and $H^1(C)$ contribute only to lower degrees.)

The motivation of this work comes from Wilson's observation
\cite{Wilson} that ``completed phase space'' for the complex
Calogero-Moser system (a cousin of the Calogero-Sutherland system) is
diffeomorphic to the Hilbert scheme of points on the affine plane.
(The author learned this observation from Segal's talk at Warwick,
1996 March.)  This made the author to look for the connection between
Hilbert schemes and Jack polynomials. However, the connection found in
this paper does not follow the Wilson's route.  Clarifying relation
between these two connections should be an important problem. For
example, it is desirable to have a geometric realization of commuting
differential operators.

It seems natural to conjecture that we get similar geometric
realization of Macdonald polynomials (two parameter family explained
above) if we replace the homology by the equivariant K-theory.
Analogous phenomenon was found for the affine Weyl group and the
affine Hecke algebra. The group ring of the affine Weyl group is
realized on the homology group of the cotangent bundle of the flag
manifold, while the Hecke algebra is realized on the equivariant
K-theory. (See \cite{Gi-book}.) We hope to return back in near
future.

\subsection*{Acknowledgment}
The author would like to thank K.~Hasegawa who told him the relation
between Jack polynomials and the Calogero-Sutherland system just after
Segal's talk at Warwick.
It is also a pleasure to acknowledge discussions with 
T.~Gocho,
A.~Matsuo and
H.~Ochiai
during the seminar on Hilbert schemes of points held in 1996 spring.
In particular, Matsuo's talk on Macdonald polynomials was very helpful.

\section{Preliminaries}\label{sec:pre}
In this section we review the theory of symmetric functions and the
equivariant cohomology for later use.

\subsection{Symmetric Functions}\label{subsec:symmetric}
First we briefly recall the theory of symmetric functions.  See
\cite{Symmetric} for detail.

A {\it partition\/} $\lambda = (\lambda_1,\lambda_2,\lambda_3,\dots)$
is a nonincreasing sequence of nonnegative integers such that
$\lambda_i = 0$ for all but finitely many $i$. Let $|\lambda| = \sum_i
\lambda_i$.  We say $\lambda$ is a partition of $n$ if $|\lambda| =
n$. We also use another presentation $\lambda = (1^{m_1}2^{m_2}\dots)$
where $m_k = \# \{i\mid \lambda_i = k\}$. Number of nonzero entries in
$\lambda$ is called {\it length\/} of $\lambda$ and denoted by
$l(\lambda)$.

If $\lambda$ and $\mu$ are partitions, we define $\lambda\cup\mu$ be
the partition whose entries are those of $\lambda$ and $\mu$ arranged
in the descending order.

For a partition $\lambda = (\lambda_1,\lambda_2,\dots)$, we give a
Young diagram such that the number of boxes in the $i$th column is
$\lambda_i$. Remark that our convention differs from one used in
\cite{Symmetric}. Our diagram is rotated by $\pi/2$ from one used in
\cite{Symmetric}.

The conjugate of a partition $\lambda$ is the partition $\lambda'$
whose diagram is the transpose of the diagram of $\lambda$, i.e.,
\begin{equation}
\label{eq:conj}
  \lambda_i' = \# \{j\mid \lambda_j\ge i\}.
\end{equation}

We define $\lambda\ge\mu$ if $|\lambda| = |\mu|$ and
\begin{equation}
  \lambda_1+\lambda_2+\cdots+\lambda_i \ge
   \mu_1+\mu_2+\cdots+\mu_i \qquad\text{for all $i$}.
\label{eq:dom}\end{equation}
This defines a partial order on the set of partitions and is called
{\it dominance order}. Note that $\lambda\ge \mu$ if and only if
$\mu'\ge \lambda'$.

Let $\Lambda_N$ be the ring of symmetric functions
\begin{equation*}
  \Lambda_N = \Z[x_1,\dots,x_N]^{{\mathfrak S}_N},
\end{equation*}
where the symmetric group ${\mathfrak S}_N$ acts by the permutation of
the variables. It is a graded ring:
\begin{equation*}
  \Lambda_N = \bigoplus_{n \ge 0} \Lambda_N^n,
\end{equation*}
where $\Lambda_N^n$ consists of the homogeneous symmetric functions
of degree $n$.
It is more relevant for us to consider symmetric functions in
``infinitely many variables'' formulated as follows: let $M\ge N$ and
consider the homomorphism
\begin{equation*}
  \Z[x_1,\dots,x_M] \to \Z[x_1,\dots,x_N]
\end{equation*}
which sends $x_{N+1}$, \dots, $x_M$ to $0$.
We have induced homomorphisms
\begin{equation*}
  \rho^n_{M,N}\colon \Lambda_M^n \to \Lambda_N^n,
\end{equation*}
which is surjective for any $M\ge N$,
and bijective for $M\ge N\ge n$.
Let
\begin{equation*}
  \Lambda^n \defeq \varprojlim \Lambda_N^n.
\end{equation*}
Then the ring of symmetric functions in infinitely many variables is 
defined by
\begin{equation*}
  \Lambda \defeq \bigoplus_n \Lambda^n.
\end{equation*}

In the relationship between symmetric functions and Hilbert schemes,
the degree $n$ corresponds the number of points, while the number of
variables $N$ are irrelevant. This is the only reason why we use the
different notation from \cite{Symmetric}.

There are several distinguished classes of symmetric functions. The
first class is the {\it monomial symmetric function\/} $m_\lambda$.
Let $\lambda$ be a partition with $l(\lambda)\le N$. Let
\begin{equation*}
  m_\lambda(x_1,\dots,x_N) 
   \defeq \sum_{\alpha\in{\mathfrak S}_N\cdot\lambda}
        x_1^{\alpha_1}\cdots x_N^{\alpha_N}
   = \frac{1}{\#\{\sigma\in{\mathfrak S}_N \mid \sigma\cdot\lambda = \lambda\}}
        \sum_{\sigma\in{\mathfrak S}_N}
        x_1^{\lambda_{\sigma(1)}}\cdots x_N^{\lambda_{\sigma(N)}},
\end{equation*}
where $\alpha = (\alpha_1,\dots,\alpha_N)$ runs over all distinct
permutation of $(\lambda_1,\lambda_2, \dots, \lambda_N)$.
If $l(\lambda) \le N$, we have
\begin{equation*}
  \rho_{M,N}^n m_\lambda(x_1,\dots,x_M) = m_\lambda(x_1,\dots,x_N).
\end{equation*}
Hence $m_\lambda$ defines an element in $\Lambda$, which is also denoted
by $m_\lambda$.
Then $\{ m_\lambda\}_{\lambda}$ is a basis for $\Lambda$.

%

The {\it $n$th power sum\/} is
\begin{equation*}
  p_n \defeq \sum x_i^n = m_{(n)}.
\end{equation*}
For a partition $\lambda = (\lambda_1,\lambda_2,\dots)$, let
$p_\lambda = p_{\lambda_1}p_{\lambda_2}\cdots$.
Then $\{ p_{\lambda} \}_{\lambda}$ is a basis for $\Lambda\otimes\Q$.
(It is {\it not\/} a $\Z$-basis for $\Lambda$.)

For a positive real number $\alpha$, we define an inner product
$\langle\cdot,\cdot\rangle$ on $\Lambda\otimes\Q$ by
\begin{equation*}
  \langle p_\lambda, p_\mu\rangle
  \defeq \alpha^{l(\lambda)} z_\lambda\delta_{\lambda\mu},
\end{equation*}
where $z_\lambda = \prod k^{m_k} m_k!$ for $\lambda =
(1^{m_1}2^{m_2}\cdots)$.

The Jack polynomials are defined by
\begin{Theorem}[\protect{\cite[10.13]{Symmetric}}]
  For each partition $\lambda$, there is a unique symmetric polynomial
  $\Jackl$ satisfying
\begin{enumerate}
\def\labelenumi{(\theequation.\theenumi)}
\def\theenumi{\arabic{enumi}}
\item\label{tri} $\Jackl = m_\lambda + \sum_{\mu < \lambda}
  u_{\lambda\mu}^{(\alpha)}m_\mu$ for suitable coefficients
  $u_{\lambda\mu}^{(\alpha)}$,
\item $\langle\Jackl, \Jackm\rangle = 0$ if $\lambda\ne\mu$.
\end{enumerate}
\label{thm:defJack}\end{Theorem}

The uniqueness is clear since the basis $\{\Jackl\}_{\lambda}$ is
obtained by the Gram-Schmidt orthogonalization from
$\{m_\lambda\}_{\lambda}$ with respect to {\it any\/} total order
compatible with the dominance order. The nontrivial point lies in 
(\ref{thm:defJack}.\ref{tri}) where the summation is over $\mu$ which
is smaller than $\lambda$ with respect to the dominance order.

Our geometric construction gives a new proof of this theorem.

The ``integral form'' $J_\lambda^{(\alpha)}$ of $\Jackl$ is defined by
the normalization
\begin{equation}
  J_\lambda^{(\alpha)} \defeq c_\lambda(\alpha)\Jackl, \qquad
  c_\lambda(\alpha) \defeq \prod_{s\in\lambda} (\alpha a(s) + l(s)+1),
\label{eq:normal}\end{equation}
where ``$s\in\lambda$'' means that $s$ is a box of the Young diagram
corresponds to $\lambda$, and $a(s)$ (resp.\ $l(s)$) is the {\it
  arm\/} (resp.\ {\it leg\/}) length defined by
\begin{equation}
\label{fig:hooklength}
\newcommand{\hf}{\hfil}
\newcommand{\hs}{\heartsuit}
\newcommand{\sps}{\spadesuit}
\Yvcentermath1
\young(\hf,\hf\hs,\hf\hs\hf,\hf s\sps\sps,\hf\hf\hf\hf)\qquad\qquad
\begin{matrix}
 a(s) &= \text{number of $\hs$} \\
 l(s) &= \text{number of $\sps$}
\end{matrix}
\end{equation}
Since our diagram is rotated, there is no reason to call them arm and
leg. But we follow the traditional convention.

The {\it augmented\/} monomial symmetric function is defined by
\begin{equation*}
  \tilde{m}_\lambda \defeq u_\lambda m_\lambda, \qquad
  u_\lambda\defeq \prod_k m_k \quad\text{for
  $\lambda = (1^{m_1}2^{m_2}\cdots)$}.
\end{equation*}

Macdonald conjectured \cite[10.26?]{Symmetric} that
$J^{(\alpha)}_\lambda$ is expressed as a linear combination of
augmented monomial symmetric function $\tilde{m}_\lambda$ with
coefficients in $\Z_{\ge 0}[\alpha]$. This conjecture was proved
affirmatively by Knop and Sahi \cite{KS} by a combinatorial method.

\subsection{Equivariant Cohomology}

We define the equivariant cohomology using the Borel construction. We
assume the group is $S^1$, though the adaptation to the general
compact Lie group is straightforward.  Let $ES^1\to BS^1$ be the
universal $S^1$-bundle which is given by the inductive limit of the
Hopf fibration $S^{2n+1}\to \CP^n$.
For a topological space $M$ with a circle action, let $M_{S^1} =
ES^1\times_{S^1} M$. We have a projection $M_{S^1} \to BS^1$, which is
a fibration with fiber $M$.  Then, the equivariant cohomology
$\HT^*(M)$ is, by definition, the cohomology of $M_{S^1}$. Similarly,
the equivariant cohomology with compact support, denoted by
$\HTc^*(M)$ is defined by the cohomology of $M_{S^1}$ with compact
support in the fiber direction of $M_{S^1}\to BS^1$.
We assume coefficients are complex numbers.

If $M$ is the space $pt$ consisting of a single point with a trivial
$S^1$-action, we have $\HT^*(pt) = H^*(BS^1)$. Since $BS^1$ is the
infinite dimensional projective space $\CP^\infty$, $H^*(BS^1)$ is the 
polynomial ring with a generator $u$ in $H^2(BS^1)$. We normalize $u$
to be the first Chern class of the tautological line bundle, i.e., the 
dual of the hyperplane bundle.

By the projection $M_{S^1}\to BS^1$, we have
$\HT^*(pt)\cong\C[u]$-module structures on $\HT^*(M)$ and $\HTc^*(M)$.

In order to relate the equivariant cohomology $\HT^*(M)$ to the
ordinary cohomology $H^*(M)$, we consider the Leray-Serre spectral
sequence associated with the fibration $M_{S^1}\to BS^1$. It is a
spectral sequence converging to $H^*(M_{S^1}) = \HT^*(M)$ with
$E_2$-term $E_2^{p,q} = H^p(BS^1)\otimes H^q(M)$. (Note $\pi_1(BS^1) =
0$.) The (decreasing) filtration is given by
\begin{equation*}
  F^p \HT^*(M) = \{ u^p\varphi \mid \varphi\in \HT^*(M) \}.
\end{equation*}
Then Kirwan proved
\begin{Theorem}\cite[5.8]{Kirwan}
Let $M$ be a compact symplectic manifold with a Hamiltonian
$S^1$-action. Then the Leray-Serre spectral sequence associated with
$M_{S^1}\to BS^1$ degenerates at the $E_2$-term. Thus we have
\begin{equation*}
  F^p \HT^{p+q}(M)/ F^{p+1}\HT^{p+q}(M) \cong H^p(BS^1)\otimes H^q(M).
\end{equation*}
\label{thm:degenerate}\end{Theorem}

Kirwan's proof works for a noncompact symplectic manifold provided
$f^{-1}((-\infty,c])$ is compact for all $c\in \R$, where $f$ is the
moment map associated with the Hamiltonian $S^1$-action
(see also \cite[Chapter~5]{Lecture}). Unfortunately, we do not have a
natural K\"ahler metric on the Hilbert scheme of points, so we could
not check the moment map satisfies the above condition.
However, we can take another route to save Kirwan's argument. In stead of
the gradient flow of the moment map, we use the $\C^*$-action which
extends the $S^1$-action. (See \cite[Chapter~7]{Lecture}.)

We have another short-cut when we assume $X$ is the total space of a
line bundle over $\CP^1$. Since the ordinary cohomology groups vanish 
in odd degree, it is obvious that the spectral sequence degenerates at 
$E_2$-term.

\section{Homology Group of the Hilbert Scheme}\label{sec:homology}

The purpose of this section is to identify the complexified ring of
symmetric functions $\Lambda\otimes\C$ with the homology group of the
Hilbert scheme of $X$. In fact, the result of this section holds when
$X$ is the total space of a line bundle $L$ over {\it any\/}
compact Riemann surface ${C}$.

Let $X$ be as above.
It is a nonsingular $2$-dimensional quasi-projective
surface containing ${C}$ as the $0$-section.


Let $\HilbX{n}$ be the Hilbert scheme parameterizing $0$-dimensional
subschemes of length $n$. By a result of Fogarty~\cite{Fog}, it is a
nonsingular $2n$-dimensional variety.

Let us consider the homology group $H_*(\HilbX{n})$ of the Hilbert
scheme. In \cite{Na-hilb}, we have constructed an action of the
Heisenberg algebra on the direct sum $\bigoplus_n H_*(\HilbX{n})$.
Let us briefly recall the construction.
For each $i\in\mathbb Z\setminus\{0\}$, let
$P_{C}[i]$ be a subvariety consisting
$(\idl_1,\idl_2)\in \coprod_n \HilbX{n-i}\times\HilbX{n}$ such that
\begin{equation}
\begin{cases}
    \cal J_1\supset \cal J_2,\;
    \text{$\Supp(\cal J_1/\cal J_2) = \{ x\}$ for some $x\in {C}$}
    & (\text{when $i > 0$}) \\
    \cal J_1\subset \cal J_2,\;
    \text{$\Supp(\cal J_2/\cal J_1) = \{ x\}$ for some $x\in {C}$}
    & (\text{when $i < 0$}) 
\end{cases},
\label{eq:corr}\end{equation}
where we consider points in $\HilbX{n}$ as ideals of $\shfO_X$.
Let $p_a$ be the projection to the $a$th factor in the product
$\HilbX{n-i}\times\HilbX{n}$. Note that the restriction of the
projection $p_2\colon P_{C}[i]\to \HilbX{n}$ is proper.
We define a homomorphism $H_*(\HilbX{n})\to H_*(\HilbX{n-i})$ by
\begin{equation*}
  \varphi\longmapsto 
    p_{1*}\left(p_2^*\varphi\cap \big[P_{C}[i]\big]\right),
\end{equation*}
where $p_2^*\varphi\cap$ means the cap product of the pull-back of the
Poincar\'e dual of $\varphi$ by $p_2^*$.
Since $p_2\colon P_{C}[i]\to \HilbX{n}$ is proper, the support of
$p_2^*\varphi\cap \big[P_{C}[i]\big]$ is compact, and hence 
$p_{1*}\left(p_2^*\varphi\cap \big[P_{C}[i]\big]\right)$ can be defined.
Moving $n$, we get an endomorphism on $\bigoplus_n H_*(\HilbX{n})$,
which we denote by the same symbol $P_{C}[i]$ for brevity.
Since $P_{C}[i]\cap (\HilbX{n-i}\times\HilbX{n})$ is a 
$(2n-i)$-dimensional subvariety
(see \cite[\S3]{Na-hilb} or \cite[\S8.3]{Lecture} for the proof),
$P_{C}[i]$ maps $H_{2n+k}(\HilbX{n})$ to
$H_{2(n-i)+k}(\HilbX{n-i})$. In particular, the middle degree part
($k=0$) is preserved. Then the main result of \cite{Na-hilb} is
the following commutator relation
\begin{equation}
  \big[P_{C}[i], P_{C}[j]\big] 
   = (-1)^{i-1}i\delta_{i+j,0}\langle {C},{C}\rangle \id.
\label{eq:comrel}\end{equation}
(The factor $(-1)^{i-1}i$, which was not determined in \cite{Na-hilb},
was given by Ellingsrud-Str\o mme~\cite{ES2}. See also
\cite[Chapter~9]{Lecture} for a proof in the spirit of this paper.)

Let $1$ be the generator of $H_0(\HilbX{0}) = \C$. Applying
$P_{C}[i]$ successively on $1$, we get a subspace in
$\bigoplus_n H_{2n}(\HilbX{n})$ which is the irreducible
representation of the Heisenberg algebra.

For each partition $\lambda=(1^{m_1}2^{m_2}\cdots)$ define
\begin{equation*}
  P^\lambda{C} \defeq
  P_{C}[-1]^{m_1} P_{C}[-2]^{m_2}\cdots 1
  \in H_{2|\lambda|}(\HilbX{|\lambda|}).
\end{equation*}
Since the representation of the Heisenberg algebra generated by
$P_{C}[i]$'s is irreducible, $P^\lambda{C}$'s are linearly
independent.

On the other hand, by the formula of G\"ottsche~\cite{Got,GS} for the
Poincar\'e polynomial $P_t(\HilbX{n})$ of $\HilbX{n}$
(see also \cite[Chapter~7]{Lecture} for the proof based on Morse
theory) we have
\begin{equation*}
  \sum_{n=0}^\infty q^n P_t(\HilbX{n}) =
  \prod_{m=1}^\infty
  \frac{(1 + t^{2m-1}q^m)^{b_1(X)}}
       {(1 - t^{2m-2}q^m)(1 - t^{2m}q^m)}\, .
\end{equation*}
Hence
\begin{equation}
  \sum_{n=0}^\infty q^n \dim H_{2n}(\HilbX{n}) =
  \prod_{m=1}^\infty \frac{1}{1 - q^m}\, .
\label{eq:char}\end{equation}
The right hand side of \eqref{eq:char} is the same as the character of
the irreducible representation of the Heisenberg algebra, hence
$P^\lambda{C}$'s span $\bigoplus_n H_{2n}(\HilbX{n})$. In other words,
$\{ P^\lambda{C}\,\}_{\lambda}$ is a basis for $\bigoplus_n
H_{2n}(\HilbX{n})$.

We identify $\bigoplus_n H_{2n}(\HilbX{n})$ with the polynomial ring
$\C[p_1,p_2,\cdots]$ (and hence with the complexified ring of
symmetric functions $\Lambda\otimes\C$) by
\begin{equation}
  p_\lambda = p_1^{m_1}p_2^{m_2}\cdots \longmapsto
  P^\lambda{C} = 
     P_{C}[-1]^{m_1} P_{C}[-2]^{m_2}\cdots 1\qquad
  \text{for $\lambda = (1^{m_1}2^{m_2}\cdots)$.}
\label{eq:ident}\end{equation}
Then the operator $P_{C}[-i]$ corresponds to the multiplication by
$p_i$ when $i > 0$.
The main result of this section is
\begin{Theorem}
  \rom{(1)} The direct sum $\bigoplus_n H_{2n}(\HilbX{n})$ of the
  middle degree homology group is isomorphic to the complexified ring
  of symmetric functions $\Lambda\otimes\C$ under the
  identification~\eqref{eq:ident}.

\rom{(2)} The intersection pairing $\langle\cdot,\cdot\rangle$ on
$H_{2n}(\HilbX{n})$ is given by
\begin{equation*}
    \langle P^\lambda{C}, P^\mu{C}\rangle
     = (-1)^n \delta_{\lambda\mu}z_\lambda
           (-\langle{C},{C}\rangle)^{l(\lambda)},
\end{equation*}
where $\lambda$, $\mu$ is a partition of $n$.
\label{thm:ident}\end{Theorem}

If we define a new inner product by
$\langle\cdot | \cdot\rangle \defeq (-1)^n\langle\cdot,\cdot\rangle$
on $H_{2n}(\HilbX{n})$, it is equal to the one used for the
definition of Jack's symmetric functions, where the parameter $\alpha$
is $-\langle{C},{C}\rangle$.

\begin{proof}[Proof of \thmref{thm:ident}]
The only remaining is to prove the statement~(2).

We identify $\bigoplus_n H_{2n}(\HilbX{n})$ with $\C[p_1,p_2,\cdots]$
by (1).
Then
the commutation relation \eqref{eq:comrel} means that the operator
$P_{C}[i]$ for $i > 0$ corresponds to
\begin{equation*}
  (-1)^{i-1} i\, \langle{C},{C}\rangle \pd{}{p_i}.
\end{equation*}
Hence for $i, j > 0$, we have
\begin{equation}
\begin{split}
  & \big[P_{C}[i]^m, P_{C}[-j]^n\big]
  = \left((-1)^{i-1} i\,\langle{C},{C}\rangle 
          \pd{}{p_i}\right)^m\; p_j^n \\
  =\;& \begin{cases}
n(n-1)\cdots (n-m+1)\delta_{ij}
      \left\{(-1)^{i-1}i\,\langle{C},{C}\rangle\right\}^m 
              P_{C}[-j]^{n-m},
&\text{for $n\ge m$}\\
0, &\text{for $n < m$.}
\end{cases}
\end{split}
\label{eq:PPcom}\end{equation}

By construction, $P_{C}[i]$ is the adjoint of $P_{C}[-i]$ with
respect to the intersection form.
For $\lambda = (1^{m_1}2^{m_2}\cdots)$, $\mu = (1^{n_1}2^{n_2}\cdots)$
let $\barlambda = (2^{m_2}3^{m_3}\cdots)$,
$\barmu = (2^{n_2}3^{n_3}\cdots)$.
Then we have
\begin{equation*}
\begin{split}
  \langle P^\lambda{C}, P^\mu{C}\rangle
=& \langle P_{C}[-1]^{m_1} P^\barlambda{C},
          P_{C}[-1]^{n_1} P^\barmu{C}\rangle \\
=& \langle P_{C}[1]^{n_1}P_{C}[-1]^{m_1}P^\barlambda{C},
          P^\barmu{C}\rangle.
\end{split}
\end{equation*}
Since $P_{C}[-1]$ commutes with $P_{C}[i]$ for $i\ne -1$, we
have
\begin{equation*}
\begin{split}
  &P_{C}[1]^{n_1}P_{C}[-1]^{m_1}P^\barlambda{C}
  = \Big[P_{C}[1]^{n_1}, P_{C}[-1]^{m_1}\Big]
    P^\barlambda{C}\\
  =\;&  \begin{cases}
  m_1!\; \langle{C},{C}\rangle^{m_1}
      P^\barlambda{C} &\text{for $m_1 = n_1$}\\
  0, &\text{otherwise,}
\end{cases}
\end{split}
\end{equation*}
where we have used \eqref{eq:PPcom}.
Hence we have
\begin{equation*}
  \langle P^\lambda{C}, P^\mu{C}\rangle = \delta_{m_1,n_1} 
  m_1! \langle{C},{C}\rangle^{m_1} 
   \langle P^\barlambda{C}, P^\barmu{C}\rangle.
\end{equation*}
Inductively, we get
\begin{equation*}
  \begin{split}
    \langle P^\lambda{C}, P^\mu{C}\rangle
 &= \delta_{\lambda\mu}
     \prod_i \left((-1)^{i-1}i
      \langle{C},{C}\rangle\right)^{m_i}\, m_i!\\
 &= \delta_{\lambda\mu} z_\lambda
     (-1)^{\sum i m_i} \left(-\langle{C},{C}\rangle\right)^{\sum m_i}.
  \end{split}
\end{equation*}
Since $\sum i m_i = n$, $\sum m_i = l(\lambda)$, we get the assertion.
\end{proof}

\section{Certain Subvarieties and Monomial Symmetric
Functions}\label{sec:subvar}

In this section we define certain middle dimensional subvarieties
parametrized by partitions and identify them with monomial symmetric
functions under the isomorphism given in \thmref{thm:ident}. The
subvarieties were first introduced by Grojnowski~\cite{Gr}, and their
identification with monomial symmetric functions was proved in
\cite[Chapter~9]{Lecture}. We reproduce them for the sake of the
reader.

Let $X$ be as in \secref{sec:homology}, i.e., the total space of a
line bundle over a compact Riemann surface $C$.  Let us consider
the $\C^*$-action on $X$ given by the multiplication on fibers. It
induces an action on the Hilbert scheme $\HilbX{n}$. Let us consider
a subvariety of the Hilbert scheme $\HilbX{n}$ defined by
\begin{equation}
  \{ Z\in\HilbX{n} \mid \Supp(\shfO_Z)\subset{C} \}.
\label{eq:LSigma}\end{equation}
Since the $\C^*$-action retracts $\HilbX{n}$ to a neighborhood of the 
above subvariety, $\HilbX{n}$ is homotopically equivalent to the subvariety.
It was pointed out by Grojnowski \cite{Gr} (see
\cite[Chapter~7]{Lecture} for the proof) that its irreducible components
$L^\lambda{C}$ are indexed by a partition $\lambda$ as follows:
first let $\pi\colon\HilbX{n}\to S^n X$ be the Hilbert-Chow morphism
which assigns to a closed subscheme $Z$ of $X$, the $0$-cycle
consisting of the points of $Z$ with multiplicities given by the
length of the local rings on $Z$. Then \eqref{eq:LSigma} is given by
$\pi^{-1}(S^n{C})$, where $S^n{C}$ is considered as a subvariety
of $S^n X$. We have a stratification of $S^n{C}$ given by
\begin{equation*}
  S^n {C} = \bigcup_{\lambda} S^n_\lambda {C}, \quad
  \text{where }S^n_\lambda {C} \defeq
    \left\{ \left.\sum_{i=1}^k \lambda_i [x_i] \in S^n {C} \right|
    \text{$x_i \neq x_j$ for $i \neq j$} \right\},
\end{equation*}
where $\lambda = (\lambda_1,\lambda_2,\dots)$ runs over
partitions of $n$.
Now consider a locally closed subvariety $\pi^{-1}(S^n_\lambda{C})$ of
\eqref{eq:LSigma}.  By a result of Brian\c{c}on~\cite{Bri}, the
Hilbert-Chow morphism $\pi$ is semismall (see also \cite{GS} or
\cite[Chapter~6]{Lecture}), the dimension of the fiber of $\pi$ over a
point in $S^n_\lambda{C}$ is $n-l(\lambda)$. Hence the dimension of
$\pi^{-1}(S^n_\lambda{C})$ is equal to $n$, which is independent of the
partition $\lambda$.  Moreover, the fiber of $\pi$ is irreducible again by
a result of Brian\c{c}on~\cite{Bri}.  Thus the irreducible components
of \eqref{eq:LSigma} are given by
\begin{equation}
   L^\lambda{C} \defeq \text{Closure of }\pi^{-1}(S^n_\lambda{C}).
\label{eq:Llambda}\end{equation}

We need another definition of $L^\lambda{C}$.  Let us consider a fixed
point $\idl$ of the $\C^*$-action in $\HilbX{n}$.  Since the fixed
point in $X$ is the zero section $C$, we have
$\Supp(\shfO_X/\idl)\subset C$.  Let us decompose $\idl$ as
$\idl_1\cap\cdots\cap\idl_m$ according to the support, i.e.,
$\Supp(\shfO_X/\idl_k) = \{x_k\}$ and $x_k\ne x_l$ for $k\ne l$. For
each $k$, we take a coordinate system $(z,\xi)$ around $x_k$ where $z$
is a coordinate of ${C}$ around $x_k$, and $\xi$ is a fiber coordinate
of $L_{x_k}$. Then $\idl_k$ is generated by monomials in $z$ and
$\xi$:
\begin{equation*}
  \idl_k = (\xi^{\lambda^{(k)}_1}, z\xi^{\lambda^{(k)}_2}, \dots,
                      z^{N-1}\xi^{\lambda^{(k)}_N}, z^N),
\end{equation*}
for some partition $\lambda^{(k)} =
(\lambda^{(k)}_1,\lambda^{(k)}_2,\dots)$ with $N = l(\lambda^{(k)})$.
For a partition $\lambda$ and a point $x\in C$, let
\begin{equation}
   \idl_{\lambda,x} \defeq (\xi^{\lambda_1}, z\xi^{\lambda_2}, \dots,
                      z^{N-1}\xi^{\lambda_N}, z^N),
\label{eq:idlLx}\end{equation}
where $N = l(\lambda)$ and $(z,\xi)$ is the coordinate system around
$x$. (See Figure~\ref{fig:Young} for the visualization of
$\idl_{\lambda,x}$.)  Thus the fixed point can be written as $\idl =
\bigcap \idl_{\lambda^{(k)},x_k}$ for some distinct points $x_k$'s in
$C$ and partitions $\lambda^{(k)}$.

If $x_k$ approaches to $x_l$, 
$\idl_{\lambda^{(k)},x_k}\cap\idl_{\lambda^{(l)},x_l}$ converges to
$\idl_{\lambda^{(k)}\cup\lambda^{(l)}, x_l}$.
This shows that the fixed point components are parametrized by
partitions $\lambda = \bigcup \lambda^{(k)}$. Let us denote by
$S^\lambda{C}$ the corresponding component. By the above discussion,
we have
\begin{equation*}
  S^\lambda{C} \defeq
  \left\{ \idl = \bigcap \idl_{\lambda^{(k)},x_k} \mid x_i\in {C},\;
  x_i \neq x_j \text{ for $i\neq j$},\;
\lambda = \bigcup \lambda^{(k)}
\right\}.
\end{equation*}
As a complex manifold, $S^\lambda{C}$
is isomorphic to a product of symmetric product
\[
S^\lambda{C}\cong S^{m_1}{C}\times S^{m_2}{C}\times\cdots,
\]
where $\lambda = (1^{m_1}2^{m_2}\dots)$.
Now another description of $L^\lambda{C}$ is
\begin{equation}
  L^\lambda{C} = \text{Closure of }
   \{ \idl\in\HilbX{n} \mid \lim_{t\to\infty} t\cdot\idl \in
   S^\lambda{C} \},
\label{eq:limit}\end{equation}
where $t\cdot$ denotes the $\C^*$-action.
In fact, since the Hilbert-Chow morphism $\pi\colon\HilbX{n}\to S^nX$
is $\C^*$-equivariant, $t\cdot\idl$ stays a compact set, or
equivalently converges to a fixed point if and only if
$\Supp(\shfO_X/\idl)\subset{C}$. Moreover,
$\lim_{t\to\infty} t\cdot\idl$ is contained in the open stratum
\[
\left\{ \idl_{(\lambda_1),x_1}\cap \dots \cap
    \idl_{(\lambda_N),x_N} \mid x_i\in {C},\;
x_i \neq x_j \text{ for $i\neq j$}\right\}
\]
of $S^\lambda{C}$
if and only if $\idl\in\pi^{-1}(S^n_\lambda{C})$. This shows the
identification \eqref{eq:limit}.

\begin{Remark}
As is presented in \cite[Chapter~7]{Lecture}, we can start from
\eqref{eq:limit} as a definition of $L^\lambda{C}$.
Studying the weight decomposition of the tangent space at a point in
$S^\lambda{C}$, we can prove that \eqref{eq:limit} is $n$-dimensional
(cf.\ discussion for \propref{prop:pullback}). The irreducibility of
\eqref{eq:limit} follows from that of $S^\lambda{C}$. Hence
\eqref{eq:limit} is an irreducible component of \eqref{eq:LSigma}.
\end{Remark}

\begin{Theorem}
  Under the isomorphism $\Lambda\otimes\C \cong\bigoplus
H_{2n}(\HilbX{n})$ given in \thmref{thm:ident},
the class $[L^\lambda{C}]$ corresponds to the monomial symmetric
function $m_\lambda$.
\label{thm:monom}\end{Theorem}

\begin{proof}
We shall show
\begin{equation}
  P_{C}[-i][L^\lambda{C}] = \sum_\mu a_{\lambda\mu} [L^\mu{C}],
  \qquad\text{for any $i \in\mathbb Z_{>0}$},
\label{eq:Ll-ind}\end{equation}
where the summation is over partitions $\mu$ of $|\lambda|+i$ which is
obtained as follows:
\begin{enumerate}
\refstepcounter{equation}\label{eq:coeff}
\def\labelenumi{(\theequation.\theenumi)}
\def\theenumi{\arabic{enumi}}
\item Add $i$ to a term in $\lambda$, say
  $\lambda_k$ \rom(possibly $0$\rom).
\item Then arrange it in descending order.
\item Define the coefficient $a_{\lambda\mu}$ by
$\#\{ l\mid \mu_l = \lambda_k + i\}$.
\end{enumerate}
This is the same as the relation between the power sum and the
monomial symmetric function:
\begin{equation*}
  p_i m_\lambda = \sum_\mu a_{\lambda\mu} m_\mu.
\end{equation*}
By induction, we have $[L^\lambda{C}] = m_\lambda$ under the
identification $P_{C}[-i] = p_i$.

The equation~\eqref{eq:Ll-ind} is proved by studying the intersection
product.
Let $p_a\colon \HilbX{n}\times\HilbX{n}\to \HilbX{n}$ be the
projection onto the $a$th factor ($a = 1,2$).
By the definition~\eqref{eq:corr}, we represent $P_{C}[-i]$ as a
subvariety. Then its set theoretical intersection with
$p_2^{-1}(L^\lambda{C})$ is
\begin{equation}
  \{ (\idl_1,\idl_2)\mid \idl_2\in L^\lambda{C},\; \idl_1\subset\idl_2,\;
  \Supp (\idl_2/\idl_1) = \{x\}
\text{ for some } x\in{C}\}.
\label{eq:intersection}\end{equation}

Let $\mu = (\mu_1,\mu_2,\dots)$ be a partition of $|\lambda|+i$
which does not necessarily satisfy the condition in \eqref{eq:coeff}.
Let $N = l(\mu)$.
Since $\{[L^\mu{C}]\}_{\mu}$ is a basis for
$H_{2|\mu|}(\HilbX{|\mu|})$, the left hand side of \eqref{eq:Ll-ind}
can be written as a linear combination of $[L^\mu{C}]$'s.
In order to determine the coefficients of
$[L^\mu{C}]$ in $P_{{C}}[-i][L^\lambda{C}]$,
it is enough to take arbitrary point $\idl_1$ in
$L^\mu{C}$ and restrict cycles to a neighborhood of $\idl_1$.
We choose the point $\idl_1 =
\idl_{(\mu_1),x_1}\cap\dots\cap\idl_{(\mu_N),x_N}\in L^\mu{C}$ where
$x_k$'s are distinct points in $C$. Here 
$\idl_{(\mu_k),x_k} = (\xi^{\mu_k}, z)$
for the coordinate $(z, \xi)$ around $x_k$.
(See \eqref{eq:idlLx} for the definition of $\idl_{\lambda,x}$.)

Suppose this point $\idl_1$ is contained in the image of
\eqref{eq:intersection} under the projection $p_1$, i.e., 
there exists $\idl_2$ such that $(\idl_1,\idl_2)$ is a point in
\eqref{eq:intersection}.
Then the point $x$ must be one of $x_k$'s, and
\begin{equation*}
   \idl_2 = \idl_{(\mu_1),x_1}\cap\dots\cap\idl_{(\mu_{k-1}),x_{k-1}}
     \cap\idl_{(\mu_k-i), x_k}
     \cap\idl_{(\mu_{k+1}),x_{k+1}}\cap\dots\cap\idl_{(\mu_N),x_N},
\end{equation*}
i.e., $\idl_2$ is obtained from $\idl_1$ by replacing
$\idl_{(\mu_k),x_k}$ by $\idl_{(\mu_k-i), x_k}$.
Since $\idl_2$ must be a point in $L^\lambda{C}$, 
$\mu$ is obtained by (a) adding $i$
to $\lambda_k$, and then (b) arranging in descending order.

Moreover, if $a_{\lambda\mu}$ is as in \eqref{eq:coeff}, there are
exactly $a_{\lambda\mu}$ choices of $x_k$'s. This explains the coefficient
$a_{\lambda\mu}$ in the formula. Thus the only remaining thing to check is 
that each choice of $(\idl_1, \idl_2)$ contributes to
$P_{C}[-i][L^\lambda{C}]$ by $[L^\mu{C}]$.
This will be shown by checking
\begin{enumerate}
\refstepcounter{equation}\label{eq:wanttoshow}
\def\labelenumi{(\theequation.\theenumi)}
\def\theenumi{\arabic{enumi}}
\item $P_{C}[-i]$ and $p_2^{-1}(L^\lambda{C})$ intersect
  transversally,
\item the intersection \eqref{eq:intersection} is isomorphic to 
$L^\mu{C}$ under the first projection $p_1$,
\end{enumerate}
in a neighborhood of $(\idl_1,\idl_2)$.

Since $\idl_1$ and $\idl_2$ are isomorphic outside $x_k$, we can
restrict our concern to $\idl_{(\mu_k),x_k}$ and $\idl_{(\mu_k-i),x_k}$.
We take the following coordinate neighborhood around 
$(\idl_{(\mu_k),x_k}, \idl_{(\mu_k-i), x_k})$ in
$\HilbX{\mu_k}\times\HilbX{\mu_k-i}$:
\begin{align*}
  &\left\{ ( (\xi^{\mu_k}+f_1(\xi), z+ g_1(\xi)),
        (\xi^{\mu_k-i} + f_2(\xi), z+g_2(\xi)))
     \mid \text{$f_1$, $g_1$, $f_2$, $g_2$ as follows}\right\} \\
  &\qquad\qquad
    f_1(\xi) = a_1\xi^{\mu_k-1} + a_2\xi^{\mu_k-2} + \dots + a_{\mu_k}, \\
  &\qquad\qquad
    g_1(\xi) = b_1 + b_2\xi + \dots + b_{\mu_k-i}\xi^{\mu_k-i-1} \\
  &\qquad\qquad\qquad\qquad\qquad
                + (b_{\mu_k-i+1} + b_{\mu_k-i+2}\xi + \dots +
                    b_{\mu_k}\xi^{i-1})(\xi^{\mu_k-i} + f_2(\xi))\\
  &\qquad\qquad
    f_2(\xi) = a'_1\xi^{\mu_k-i-1} 
                 + a'_2\xi^{\mu_k-i-2} + \dots + a'_{\mu_k-i}, \\
  &\qquad\qquad
    g_2(\xi) = b'_1+ b'_2\xi + \dots + b'_{\mu_k-i}\xi^{\mu_k-i-1}
\end{align*}
where $(a_1,\dots, a_{\mu_k}, b_1,\dots, b_{\mu_k})$ 
(resp.\ $(a'_1,\dots, a'_{\mu_k-i}, b'_1,\dots, b'_{\mu_k-i})$) 
is in a neighborhood of $0$ in $\C^{2\mu_k}$ (resp.\ $\C^{2(\mu_k-i)}$).
Then the above ideal is contained in
$P_{C}[-i]$ if and only if the followings hold
\begin{enumerate}
\refstepcounter{equation}
\def\labelenumi{(\theequation.\theenumi)}
\def\theenumi{\arabic{enumi}}
\item $\xi^{\mu_k} + f_1(\xi) = \xi^i(\xi^{\mu_k-i} + f_2(\xi))$,
\item $g_1(\xi) - g_2(\xi)$ is divisible by $\xi^{\mu_k-i} + f_2(\xi)$.
\end{enumerate}
Namely, the defining equation for $P_{{C}}[-i]$ is
\begin{equation*}
\begin{split}
  & a_1 = a'_1, a_2 = a'_2, \dots, a_{\mu_k-i} = a'_{\mu_k-i}, \\
  & a_{\mu_k-i+1} = \cdots = a_{\mu_k} = 0, \\
  & b_1 = b'_1, b_2 = b'_2, \dots, b_{\mu_k-i} = b'_{\mu_k-i}.
\end{split}
\end{equation*}

On the other hand, the defining equation for
$p_2^{-1}(L^\lambda{C})$ is
\begin{equation*}
  a'_1 = a'_2 = \cdots = a'_{\mu_k-i} = 0.
\end{equation*}
Now our assertions \eqref{eq:wanttoshow} are immediate.
\end{proof}

Our next task is to explain a geometric meaning of the dominance
order \eqref{eq:dom}. It is given by modifying the stratification
introduced in \cite{Bri,Iar}.

For $i\ge 0$, let $\shfO_X(-iC)$ be the sheaf of functions vanishing
to order $\ge i$ along the zero section $C$. Let $\idl\subset\shfO_X$
be an ideal of colength $n$ such that the support of $\shfO_X/\idl$ is
contained in $C$.  We consider the sequence
$(\lambda'_1,\lambda'_2,\dots)$ of nonnegative integers given by
\begin{equation*}
  \lambda'_i(\idl)\defeq \operatorname{length}\left(
   \frac{\shfO_X(-(i-1)C)}{\idl\cap\shfO_X(-(i-1)C)+\shfO_X(-iC)}
  \right).
\end{equation*}
The reason why we put the prime become clear later.  The sequence in
\cite{Bri,Iar} was defined by replacing $\shfO_X(-iC)$ by 
$\frak m_x^i$ where $\frak m_x$ is the maximal ideal corresponding to
a point $x$.  As in \cite[Lemma~1.1]{Iar}, we have
$\idl\supset\shfO_X(-nC)$, hence $\lambda'_i(\idl) = 0$ for $i\ge
n+1$.  From the exact sequence
\begin{equation*}
 0 \to \frac{\shfO_X(-iC)}{\idl\cap\shfO_X(-iC)} \to
      \frac{\shfO_X(-(i-1)C)}{\idl\cap\shfO_X(-(i-1)C)} \to
 \frac{\shfO_X(-(i-1)C)}{\idl\cap\shfO_X(-(i-1)C)+\shfO_X(-iC)}
 \to 0,
\end{equation*}
we have
\begin{equation*}
  \sum_{i=1}^n \lambda'_i(\idl) = n.
\end{equation*}

Let us decompose the ideal $\idl$ by its support, i.e.,
$\idl = \idl_1\cap\dots\cap\idl_N$ such that
$\{ \Supp(\shfO_X/\idl_k) \}_k$ are $N$ distinct points.
By definition, we have
\begin{equation}
  \lambda_i'(\idl) = \sum_{k=1}^N \lambda_i'(\idl_k).
\label{eq:sum}\end{equation}

Suppose that $\idl$ satisfies $\Supp(\shfO_X/\idl) = \{x\}$ for some
$x\in C$. We take a coordinate system $(z,\xi)$ around $x$ where $\xi$
is a coordinate for the fiber. 
If $\xi^i f_1(z),\dots, \xi^i f_d(z)$ form a basis of
\begin{equation*}
\frac{\shfO_X(-iC)}{\idl\cap\shfO_X(-iC)+\shfO_X(-(i+1)C)},
\end{equation*}
Then $\xi^{i-1} f_1(z),\dots, \xi^{i-1} f_d(z)$ are linearly
independent in 
\begin{equation*}
\frac{\shfO_X(-(i-1)C)}{\idl\cap\shfO_X(-(i-1)C)+\shfO_X(-iC)}.
\end{equation*}
Hence we have $\lambda'_i(\idl) \ge \lambda'_{i+1}(\idl)$. Thus
$(\lambda'_1(\idl), \lambda'_2(\idl), \dots)$ is a partition of $n$.
By \eqref{eq:sum}, the same is true for general $\idl$ which do not
necessarily satisfy $\Supp(\shfO_X/\idl)=\{x\}$.
Let us denote the partition by $\lambda'(\idl)$.

For a partition $\lambda' = (\lambda'_1,\lambda'_2,\dots)$ of $n$,
let $W^{\lambda'}$ be the set of ideals $\idl\subset\shfO_X$ with
colength $n = |\lambda'|$ such that $\shfO_X/\idl$ is supported on
$C$ and $\lambda'(\idl) = \lambda'$.
Since
\begin{equation*}
  \operatorname{length}(\shfO_X(-iC)/\idl\cap\shfO_X(-iC)) \leq
  \sum_{j=i+1}^n \lambda'_j = n - \sum_{j=1}^i \lambda'_j
\end{equation*}
is a closed condition on $\idl$, the union
\begin{equation*}
   \bigcup_{\mu'\ge \lambda'} W^{\mu'}
\end{equation*}
is a closed subset of $\{\idl\in \HilbX{n}\mid
\Supp(\shfO_X/\idl)\subset C\}$. Thus we have
\begin{equation}
\label{eq:closure}
  \text{Closure of }W^{\lambda'} \subset\bigcup_{\mu'\ge \lambda'}
  W^{\mu'}.
\end{equation}

Suppose that $\lambda'$ is the conjugate of $\lambda$ as in
\eqref{eq:conj}.
We get the following third description of $L^\lambda{C}$.
\begin{Proposition}
  $L^\lambda{C} = \text{Closure of }W^{\lambda'}$.
\label{prop:third}\end{Proposition}

\begin{proof}
Let us write $\lambda = (\lambda_1,\dots,\lambda_N)$ with 
$N = l(\lambda)$.
By \eqref{eq:limit}, a generic point $\idl$ in $L^\lambda{C}$ satisfies
$\lim_{t\to\infty} t\cdot\idl = 
\idl_{(\lambda_1),x_1}\cap \dots \cap \idl_{(\lambda_N),x_N}$
such that $x_i$'s are distinct points in $C$.
Since the support of $\idl$ cannot move as $t\to\infty$, we can
decompose
$\idl = \idl_1\cap\dots\cap\idl_N$ such that
$\Supp(\shfO_X/\idl_k) = \{x_k\}$.

Take a coordinate system $(z,\xi)$ around $x_k$ as before. Then
$\idl_{(\lambda_k),x_k}$ was defined by $(z,\xi^{\lambda_k})$.
Since $\lim_{t\to\infty} t\cdot\idl_k = \idl_{(\lambda_k),x_k}$, we have
\begin{equation*}
  \idl_k = (\xi^{\lambda_k}, z + a_1 \xi + a_2 \xi^2 + \cdots +
                  a_{\lambda_k-1}\xi^{\lambda_k-1})
\end{equation*}
for some $a_1,\dots, a_{\lambda_k-1}$
(cf.\ the proof for \thmref{thm:monom}).
Then $\lambda'(\idl_k) = (1^{\lambda_k})$. By \eqref{eq:sum}, we have
$$
  \lambda'(\idl) = \lambda'.
$$
This shows $L^{\lambda}{C}\subset\text{Closure of }W^{\lambda'}$.

Conversely, suppose $\idl$ is a point in $W^{\lambda'}$. Let
\begin{equation*}
  \idl_i \defeq
  \frac{\idl\cap\shfO_X(-(i-1)C)}{\idl\cap\shfO_X(-iC)}.
\end{equation*}
Then $\idl' \defeq \bigoplus \idl_i$ satisfies $\lim_{t\to\infty}
t\cdot\idl = \idl'$ and $\idl'\in W^{\lambda'}$. This shows
$W^{\lambda'}\subset L^\lambda{C}$.
\end{proof}

\section{Equivariant Cohomology of Hilbert Schemes}

Let $X$ be the total space of a line bundle over $\CP^1$ and $C$ the
zero section. $X$ is the quotient space of $(\C^2\setminus\{0\})\times\C$
by the $\C^*$-action given by
\begin{equation*}
  (z_0, z_1, \xi) \mapsto (\lambda z_0, \lambda z_1, \lambda^{-\alpha}\xi)
  \qquad \lambda\in\C^*,
\end{equation*}
where $\alpha = -c_1(L)[{C}] = -\langle{C},{C}\rangle$.
We denote by $[(z_0, z_1, \xi)]$ the equivalence class containing
$(z_0, z_1, \xi)$. The projection $X \to \CP^1$ is given by
$[(z_0, z_1, \xi)]\mapsto [z_0:z_1]$, hence $\xi$ is a coordinate
for the fiber.

We consider the $S^1$-action on $\HilbX{n}$ induced by the action on
$X$ defined by
\begin{equation}
  [(z_0, z_1, \xi)] \mapsto [(z_0, t^{-1} z_1, t^{\alpha}\xi)]\qquad
  \text{for $t\in S^1$}.
\label{eq:action}\end{equation}
Note that this differs from the $S^1$-action studied in the previous
section given by 
$[(z_0, z_1, \xi)] \mapsto [(z_0, z_1, t\xi)]$.
We shall study the equivariant cohomology of $\HilbX{n}$ with respect
to the above $S^1$-action.
The reason for using this $S^1$-action is to identify the
normalization factor~\eqref{eq:normal} with the equivariant Euler
class (see \propref{prop:pullback}).

Let $pt$ be the space consisting of a single point with a trivial
$S^1$-action. Let us denote the obvious $S^1$-equivariant morphism by
$p\colon \HilbX{n}\to pt$.
We have an equivariant push-forward $p_*$.
Define a bilinear form
$\langle \cdot, \cdot\rangle_{S^1}\colon
\HTc^k(\HilbX{n})\otimes \HTc^k(\HilbX{n})\to\HT^{2(k-2n)}(pt)$ by
\begin{equation}
  \langle \varphi,\psi\rangle_{S^1} \defeq p_*(\varphi\cup \psi),
\label{eq:inner}\end{equation}
where $\cup$ is the cup product
\begin{equation*}
  \cup\colon\HTc^k(\HilbX{n})\otimes\HTc^k(\HilbX{n})\to
     \HTc^{2k}(\HilbX{n}).
\end{equation*}
This is trivial unless $k$ is even and $k\ge 2n$. We assume this
condition henceforth.
In particular, $\langle\cdot,\cdot\rangle_{S^1}$ is
symmetric.

Let $j^*\colon \HTc^k(\HilbX{n})\to \HT^k(\HilbX{n})$ be the natural
homomorphism.
Then the cup product $\cup$ factor through $j^*$:
\begin{equation*}
  \varphi\cup \psi = \varphi\cup j^*\psi,
\end{equation*}
where $\cup$ in the right hand side is the cup product
\begin{equation*}
  \cup \colon \HTc^k(\HilbX{n})\otimes\HT^k(\HilbX{n})\to
   \HTc^{2k}(\HilbX{n}).
\end{equation*}
In particular, $\langle\cdot,\cdot\rangle_{S^1}$ is well-defined on
\begin{equation*}
  \HTc^k(\HilbX{n})/\Ker j^* \cong \Ima j^*.
\end{equation*}

Let $\HilbX{n}_{S^1} \to BS^1$ be the fibration used in the definition
of the equivariant cohomology.  By \thmref{thm:degenerate} (see also
the discussion after \thmref{thm:degenerate}), the spectral sequence of
the fibration in cohomology degenerates.  Let us denote by
$F^*\HT^k(\HilbX{n})$ the corresponding decreasing filtration.

Since 
\(H^q(\HilbX{n}) = 0\) for $q > 2n$,
we have
\begin{equation*}
  \begin{split}
  & F^p \HT^k(\HilbX{n}) = \HT^k(\HilbX{n}) \qquad\text{for $p\le k-2n$}, \\
  & \HT^k(\HilbX{n})/F^{k-2n+1}\HT^k(\HilbX{n}) \cong
    \HT^{k-2n}(pt)\otimes H^{2n}(\HilbX{n}).
  \end{split}
\end{equation*}

If $j^*\psi \in F^{k-2n+1}\HT^k(\HilbX{n})$, then
$\varphi\cup j^* \psi\in F^{2(k-2n)+1}\HTc^{2k}(\HilbX{n})$.
Hence
$p_*(\varphi\cup j^*\psi)\in F^{2(k-2n)+1}\HT^{2(k-2n)}(pt)\cong 0$.
In particular, $\langle\cdot,\cdot\rangle_{S^1}$
is well-defined on
$\Ima j^*/\Ima j^*\cap F^{k-2n+1}\HT^k(\HilbX{n})$.
Consider the composition
\begin{equation}
\begin{gathered}
  \HTc^{k}(\HilbX{n})/\operatorname{Rad}\langle\ ,\ \rangle_{S^1}
  \xrightarrow{j^*} \Ima j^*/\Ima j^*\cap F^{k-2n+1}\HT^k(\HilbX{n})\\
  \to 
  \HT^k(\HilbX{n})/F^{k-2n+1}\HT^k(\HilbX{n}) \cong
  \HT^{k-2n}(pt)\otimes H^{2n}(\HilbX{n}),
\end{gathered}
\label{eq:compo}\end{equation}
which is injective by the above discussion.

We have an analogous injective homomorphism between {\it ordinary\/}
cohomology groups
\begin{equation*}
  H^{2n}_c(\HilbX{n})/\operatorname{Rad}\langle\ ,\ \rangle
   \to H^{2n}(\HilbX{n}).
\end{equation*}
However, the bilinear form $\langle\ ,\ \rangle$ on
$H^{2n}_c(\HilbX{n})$ is nondegenerate by \thmref{thm:ident}(2).
Hence $\operatorname{Rad}\langle\ ,\ \rangle = 0$. Moreover, the
Poincar\'e duality implies
$\dim H^{2n}_c(\HilbX{n}) = \dim H^{2n}(\HilbX{n})$.
Therefore the natural homomorphism
$ H^{2n}_c(\HilbX{n})\to H^{2n}(\HilbX{n})$ is an isomorphism.

\begin{Theorem}
Suppose $k$ is even and $k\ge 2n$. Then the composition~\eqref{eq:compo}
\begin{equation}
  \HTc^{k}(\HilbX{n})/\operatorname{Rad}\langle\ ,\ \rangle_{S^1} \to
    \HT^{k-2n}(pt)\otimes H^{2n}(\HilbX{n})
\label{eq:isom}\end{equation}
is an isomorphism.
Moreover, the bilinear form $\langle\cdot ,\cdot \rangle_{S^1}$ on the left
hand side is equal to one on the right hand side 
induced by the intersection pairing on $H^{2n}(\HilbX{n})\cong
H^{2n}_c(\HilbX{n})$.
\label{thm:isom}\end{Theorem}

\begin{proof}
The subvariety $L^\lambda{C}$ defined in \eqref{eq:Llambda}
is invariant under the $S^1$-action.
Hence its Poincar\'e dual defines an element in the equivariant cohomology
$\HTc^{2n}(\HilbX{n})$, which we denote by $[L^\lambda{C}]_{S^1}$.
Under \eqref{eq:compo} with $k=2n$, 
$[L^\lambda{C}]_{S^1}\bmod\operatorname{Rad}\langle\ ,\ \rangle_{S^1}$
is mapped to $1\otimes [L^\lambda{C}]$.
More generally, if $u$ is the generator of $\HT^*(pt)$ which lives in
$\HT^{2}(pt)$, $u^{k/2 - n}[L^\lambda{C}]_{S^1}$ is mapped to
$u^{k/2 - n}\otimes [L^\lambda{C}]$.
This shows that \eqref{eq:isom} is surjective.
As we have already seen the injectivity, it is an isomorphism.

Next consider the bilinear forms on the both hand sides of
\eqref{eq:isom}. The doubt arises from that 
$\HT^k(\HilbX{n})\cong \bigoplus_{p+q=k}\HT^p(pt)\otimes
H^q(\HilbX{n})$ is not an {\it algebra\/} isomorphism, it is only an
isomorphism between {\it vector spaces}.
If we consider the graded space $G\HT^k(\HilbX{n})$ instead of
$\HT^k(\HilbX{n})$, it becomes an algebra isomorphism.
However, as we pointed out above, $F^{k-2n+1}\HT^k(\HilbX{n})$ does
not contribute to $\langle\cdot ,\cdot \rangle_{S^1}$. Hence we have
the assertion.
\end{proof}

\section{Localization}

We assume $X$ is the total space of a line bundle over $\CP^1$ with
$\alpha = -c_1(L)[\CP^1] > 0$.
Let us consider the $S^1$-action given in \eqref{eq:action}. The fixed
point set consists of two components
\begin{equation*}
  \{ [(1,0,0)] \} \text{ and } \{ [(0, 1,\xi)] \mid \xi\in\C \}.
\end{equation*}
The first is an isolated point and the latter is isomorphic to $\C$.
If $Z\in \HilbX{n}$ is fixed by the induced action, its support must
be contained in the above fixed point set.

\begin{figure}[htbp]
\begin{center}
\leavevmode
\setlength{\unitlength}{0.01250000in}
\begin{picture}(73,114)(0,-10)
\path(22,62)(22,82)(2,82)
        (2,62)(22,62)
\path(22,42)(22,62)(2,62)
        (2,42)(22,42)
\path(42,42)(42,62)(22,62)
        (22,42)(42,42)
\path(22,22)(22,42)(2,42)
        (2,22)(22,22)
\path(42,22)(42,42)(22,42)
        (22,22)(42,22)
\path(22,2)(22,22)(2,22)
        (2,2)(22,2)
\path(42,2)(42,22)(22,22)
        (22,2)(42,2)
\path(62,2)(62,22)(42,22)
        (42,2)(62,2)
\put(7,87){\makebox(0,0)[lb]{\smash{$\xi^4$}}}
\put(25,67){\makebox(0,0)[lb]{\smash{$z\xi^3$}}}
\put(45,27){\makebox(0,0)[lb]{\smash{$z^2\xi$}}}
\put(67,7){\makebox(0,0)[lb]{\smash{$z^3$}}}
\end{picture}
\caption{Young diagram and an ideal}
\label{fig:Young}
\end{center}
\end{figure}
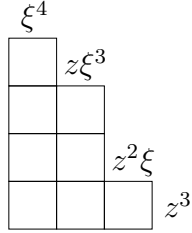

There are a family of fixed points $Z_\lambda$ with support $[(1,0,0)]$
indexed by a partition
$\lambda = (\lambda_1,\lambda_2,\dots)$ of $n$ as follows: let
$(z,\xi)$ be a coordinate system around $[(1,0,0)]$ given by $(z,\xi)
\leftrightarrow [(1,z,\xi)]$. 
Then the corresponding ideal $\idl_\lambda$ is given by
\begin{equation*}
  \idl_\lambda \defeq (\xi^{\lambda_1}, z\xi^{\lambda_2}, \dots,
                       z^{N-1}\xi^{\lambda_N}, z^N),
\end{equation*}
where $N = l(\lambda)$. Namely, $\idl_\lambda =
\idl_{\lambda,[(1,0,0)]}$ in the notation~\eqref{eq:idlLx}. The ideal
can be visualized by the Young diagram as follows. (Remind that our
diagram is rotated by $\pi/2$ from one used in \cite{Symmetric}.) If
we write the monomial $z^{p-1}\xi^{q-1}$ at the position $(p,q)$, the
ideal is generated by monomials outside the Young diagram. For
example, the Young diagram in figure~\ref{fig:Young} corresponds to
the ideal $(\xi^4, z\xi^3, z^2\xi, z^3)$. Although it is not necessary
for our purpose, one can show that these ideals are a complete list of
fixed points in $\HilbX{n}$ whose support are $[(1,0,0)]$ when $\alpha
> 0$.

Let $i_\lambda\colon Z_\lambda\to \HilbX{n}$ be the inclusion.
We have a homomorphism between equivariant cohomology groups:
\begin{equation*}
  \HT^*(Z_\lambda) \ni \varphi
   \longmapsto i_{\lambda*}\varphi \in \HTc^*(\HilbX{n}).
\end{equation*}
Let $e(T_{Z_\lambda}\HilbX{n})$ be the $S^1$-equivariant Euler class
of the tangent space at $Z_\lambda$, considered as an element of
$\HT^{4n}(Z_\lambda)$. Then we have
\begin{Theorem}
Let $\langle\cdot,\cdot\rangle_{S^1}$ be the bilinear pairing defined in
\eqref{eq:inner}. Then we have
\begin{equation*}
   \langle i_{\lambda*}\varphi, i_{\mu*}\psi\rangle_{S^1}
   = 
   \begin{cases}
    0 & \text{if $\lambda\ne\mu$}, \\
    p_{\lambda*} \left(\varphi\cup \psi\cup 
      e(T_{Z_\lambda}\HilbX{n})\right) &\text{if $\lambda=\mu$},
   \end{cases}
\end{equation*}
where $p_\lambda\colon Z_\lambda\to pt$ is the obvious map.
\label{thm:push}\end{Theorem}

\begin{proof}
By the projection formula, we have
\begin{equation*}
  i_{\lambda*}\varphi\cup i_{\mu*}\psi
 = i_{\mu*}\left( i_{\mu}^*i_{\lambda*}\varphi \cup \psi\right).
\end{equation*}
The homomorphism $i_{\mu}^*i_{\lambda*}$ is zero when
$\lambda\neq\mu$, and is given by the multiplication by the
equivariant Euler class $e(T_{Z_\lambda}\HilbX{n})$ if $\lambda=\mu$.
Therefore,
\makeatletter\tagsleft@false\begin{align}
   & \langle i_{\lambda*}\varphi, i_{\lambda*}\psi\rangle_{S^1}
  = p_*\left(i_{\lambda*}\varphi\cup i_{\lambda*}\psi\right)
  = p_*i_{\lambda*}\left(\varphi \cup \psi \cup 
      e(T_{Z_\lambda}\HilbX{n})\right) \notag \\
  =\;& p_{\lambda*}\left(\varphi \cup \psi \cup 
      e(T_{Z_\lambda}\HilbX{n})\right)
  \tag*{\qed}
\end{align}\makeatother
\renewcommand{\qed}{}
\end{proof}

Taking a class from $\HT^*(Z_\lambda)$ for each $\lambda$, we can give
a set of mutually orthogonal elements in $\HTc^*(\HilbX{n})$. These
are candidates for Jack polynomials.  Our remaining task is to choose
classes so that two conditions in \thmref{thm:defJack} will be
satisfied. Thus we need to know monomial symmetric functions, or their
geometric counterparts $[L^\lambda{C}]_{S^1}$.

In order to study $[L^\lambda{C}]_{S^1}$, we need to know the
$S^1$-module structure of its normal bundle at $Z_\lambda$.
In fact, the formula will become clearer if we
consider the two dimensional torus action given by
\begin{equation*}
    [(z_0, z_1, \xi)] \mapsto [(z_0, t_1 z_1, t_2\xi)]\qquad
  \text{for $(t_1, t_2)\in T = S^1\times S^1$}.
\end{equation*}
The $S^1$-action studied in \secref{sec:subvar} is given by the
restriction to $(1,t)\in T$, while the $S^1$-action given by
\eqref{eq:action} is the restriction to $(t^{-1}, t^{\alpha})$.
The fixed point $Z_\lambda$ defined above is also fixed by
this torus action.
The tangent space $T_{Z_\lambda}\HilbX{n}$ at $Z_\lambda$ is a
$T$-module and have a weight decomposition
\begin{equation*}
   T_{Z_\lambda}\HilbX{n} = \bigoplus_{p,q\in\mathbb Z} H(p,q),
\end{equation*}
where $H(p,q) = \{ v\in T_{Z_\lambda}\HilbX{n} \mid
(t_1,t_2)\cdot v = t_1^p t_2^q v\;\text{for $(t_1,t_2)\in T$}\}$.

\begin{Lemma}
The character of the tangent space $T_{Z_\lambda}\HilbX{n}$ is given
by
\begin{equation*}
 \sum_{p,q} \dim_\C H(p,q) T_1^p T_2^q =
 \sum_{s\in\lambda} \left(T_1^{l(s)+1}T_2^{-a(s)}
     + T_1^{-l(s)}T_2^{a(s)+1}\right).
\end{equation*}
\label{lem:weight}\end{Lemma}

\begin{proof}
The corresponding formula was proved in \cite{ES} (see also
\cite[2.2.5]{Go-book} or \cite[Proposition~5.5]{Lecture})
when $X$ is replaced by $\C^2$ with the torus action given by
\begin{equation*}
  (z_1,z_2) \mapsto (t_1 z_1, t_2 z_2)\qquad
  \text{for $(z_1,z_2)\in\C^2$, $(t_1,t_2)\in T$.}
\end{equation*}
(The presentation in \cite{ES} or \cite{Go-book} is different from the
above. After substituting the arm/leg-length to their formula, we get the
above formula.)
Since the exponential map gives a 
$T$-equivariant isomorphism between a neighborhood of $0\in\C^2$
and a neighborhood of $[(1,0,0)]\in X$, we get the assertion.
\end{proof}

Since the $S^1$-equivariant Euler class of an $S^1$-module is given by
the product of weights, we get the following.
\begin{Corollary}
The $S^1$-equivariant Euler class of the tangent space at $Z_\lambda$
is given by
\begin{equation*}
  e(T_{Z_\lambda}\HilbX{n}) =
  u^{2n}\prod_{s\in\lambda} \left(-\alpha a(s) - l(s) - 1 \right)\cdot
        \prod_{s\in\lambda} \left(\alpha(a(s) + 1) + l(s)\right).
\end{equation*}
\label{cor:tangent}\end{Corollary}

In particular, $Z_\lambda$ is an isolated fixed point of $S^1$-action
and $e(T_{Z_\lambda}\HilbX{n})$ is nonzero since $\alpha > 0$.
In particular, we have
\begin{Corollary}
Let $\varphi_\lambda$ be a nonzero element in $\HT^k(Z_\lambda)$. Then
$\{ i_{\lambda*}\varphi_\lambda\pmod{\operatorname{Rad}\langle\ ,\ 
\rangle} \}_{\lambda}$ forms an orthogonal basis of
$\HTc^{4n+k}(\HilbX{n})/\operatorname{Rad}\langle\ ,\ \rangle$.
\label{cor:orth}\end{Corollary}

Let us decompose the tangent space $T_{Z_\lambda}\HilbX{n}$ as
\begin{gather*}
  N_{Z_\lambda}^{>0} \defeq \bigoplus_{p\in\mathbb Z}
                         \bigoplus_{q > 0}H(p,q), \qquad
  N_{Z_\lambda}^{\le 0} \defeq \bigoplus_{p\in\mathbb Z}
                         \bigoplus_{q \le 0}H(p,q), \\
  T_{Z_\lambda}\HilbX{n} = N_{Z_\lambda}^{>0}\oplus N_{Z_\lambda}^{\le
    0}
\end{gather*}

\begin{Proposition}
The pull-back of the equivariant cohomology class
$[L^\lambda{C}]_{S^1}$ by the inclusion $i_{\lambda}\colon
Z_\lambda\to \HilbX{n}$ is equal to the $S^1$-equivariant Euler
class of $N_{Z_\lambda}^{>0}$. It is given by
\begin{equation*}
  i_\lambda^*[L^\lambda{C}]_{S^1} = e(N_{Z_\lambda}^{>0}) =
  u^{n}\prod_{s\in\lambda} \left(\alpha(a(s) + 1) + l(s)\right).
\end{equation*}
\label{prop:pullback}\end{Proposition}

\begin{proof}
By the description~\eqref{eq:limit} of $L^\lambda{C}$,
it is nonsingular at $Z_\lambda$ and its
tangent space is equal to $N_{Z_\lambda}^{\le 0}$.
Hence the normal bundle is $N_{Z_\lambda}^{>0}$.
Thus we have $i_\lambda^*[L^\lambda{C}]_{S^1} =
e(N_{Z_\lambda}^{>0})$. The above formula follows from
\lemref{lem:weight}.
\end{proof}

\begin{Theorem}
For each partition $\lambda$, let $F_\lambda$ be the cohomology class
\begin{equation}
i_{\lambda*}\left(\frac{u^n}{e(N_{Z_\lambda}^{\le 0})}\right) =
i_{\lambda*}\left( \frac{1}{
    \prod_{s\in\lambda} \left(-\alpha a(s) - l(s) - 1 \right)}\right)
\label{eq:Flambda}\end{equation}
considered as an element in
\begin{equation*}
  \HTc^{4n}(\HilbX{n})/\operatorname{Rad}\langle\ ,\ \rangle \cong
      \HT^{2n}(pt)\otimes H^{2n}_c(\HilbX{n})
  \xrightarrow[u^{-n}\otimes\text{P.D.}]{\cong}
  H_{2n}(\HilbX{n}).
\end{equation*}
Then they satisfy
\begin{align}
  & F_{\lambda} = [L^\lambda{C}] 
             + \sum_{\mu < \lambda} u_{\lambda\mu}^{(\alpha)}[L^\mu{C}]
    \qquad\text{for some $u_{\lambda\mu}^{(\alpha)}$,} \label{eq:FL}\\
  &\langle F_{\lambda}, F_{\mu}\rangle = 0
    \qquad\text{for $\lambda\neq\mu$}. \label{eq:orth}
\end{align}
\label{thm:Jack}\end{Theorem}

\begin{Corollary}
Under the identification $\bigoplus H_{2n}(\HilbX{n})\cong\Lambda\otimes\C$
given in \thmref{thm:ident}, $F_{\lambda}$ corresponds to the Jack's
symmetric function $P_{\lambda}^{(\alpha)}$, where the parameter
$\alpha$ is given by $-\langle{C},{C}\rangle$ \rom(see
\eqref{eq:selfinter}\rom).
\end{Corollary}

\begin{proof}[Proof of \protect{\thmref{thm:Jack}}]
The equation \eqref{eq:orth} follows from \thmref{thm:push}.

Let us take the generator $\varphi_\lambda\in\HT^0(Z_\lambda)$.  By
\corref{cor:orth}, we can write $[L^\lambda{C}]_{S^1}$ as a linear
combination of $i_{\mu*}\varphi_\lambda$ modulo
$\operatorname{Rad}\langle\ ,\ \rangle$.  One can check that $Z_\mu$
is contained in $W^{\mu'}$.  Hence $L^\lambda{C}=\text{Closure of }
W^{\lambda'}$ (see \propref{prop:third}) contains $Z_\mu$ only if
$\mu\le\lambda$ by \eqref{eq:closure}. In particular, the pull-back
$i_\mu^*[L^\lambda{C}]_{S^1} = 0$ unless $\mu\le\lambda$. Hence we can
write
\begin{equation*}
  [L^\lambda{C}]_{S^1} = \sum_{\mu\le\lambda} v_{\lambda\mu}^{(\alpha)}
  i_{\mu*}\varphi_\lambda\mod{\operatorname{Rad}\langle\ ,\ \rangle}
\end{equation*}
for some $v_{\lambda\mu}^{(\alpha)}$.

In order to compute $v_{\lambda\lambda}^{(\alpha)}$, we consider the 
pull-back by $i_\lambda^*$. We have
\begin{align*}
  i_\lambda^*[L^\lambda{C}]_{S^1} &= e(N_{Z_\lambda}^{>0})
    &&\text{(by \propref{prop:pullback})} \\
  &= \frac{e(T_{Z_\lambda}\HilbX{n})}
      {e(N_{Z_\lambda}^{\le 0})}
  &&\text{(since $T_{Z_\lambda}\HilbX{n} =
             N_{Z_\lambda}^{>0}\oplus N_{Z_\lambda}^{\le 0}$)} \\
  &=\frac{1}{u^n}i_\lambda^* i_{\lambda*}\left(
   \frac{u^n}{e(N_{Z_\lambda}^{\le 0})}\right)
  &&\text{(by the self-intersection formula)}.
\end{align*}
Thus we get the assertion.
\end{proof}

Finally note that the denominator in \eqref{eq:Flambda} coincides with
the normalization factor~\eqref{eq:normal} up to sign. Hence the
integral form $J_\lambda^{(\alpha)}$ is an {\it integral\/} class.
This shows that $J_\lambda^{(\alpha)}$ is expressed as a linear
combination of $m_\lambda$ with coefficients in $\Z$ for any positive
integer $\alpha$. Thus we have $J_\lambda^{(\alpha)} =
\sum_{\mu<\lambda} \tilde{u}_{\lambda\mu}^{(\alpha)} m_\mu$ with
$\tilde{u}_{\lambda\mu}^{(\alpha)}\in\Z[\alpha]$. This proves a part
of Macdonald conjecture \cite[10.26?]{Symmetric} mentioned in
\subsecref{subsec:symmetric}.

\end{document}